\documentclass{amsart}
\usepackage{amscd,amssymb,hyperref}

\theoremstyle{plain}
\newtheorem{thm}[subsection]{Theorem}
\newtheorem{prop}[subsection]{Proposition}
\newtheorem{lem}[subsection]{Lemma}
\newtheorem{cor}[subsection]{Corollary}

\theoremstyle{remark}
\newtheorem{rem}[subsection]{Remark}
\newtheorem{exm}[subsection]{Example}

\numberwithin{equation}{section}

\renewcommand{\b}[1]{\mathbf{#1}}

\newcommand{\A}{\mathcal{A}}
\newcommand{\LL}{\mathcal{L}}
\newcommand{\Ai}{{\mathcal A}_\infty}

\newcommand{\la}{\lambda}
\newcommand{\bul}{\bullet}
\newcommand{\C}{\mathbb{C}}
\newcommand{\Z}{\mathbb{Z}}
\newcommand{\T}{(\mathbb{C}^*)^N}
\newcommand{\CP}{{\mathbb{CP}}}
\newcommand{\ii}{{\mathrm{\,i\,}}}
\newcommand{\SL}{{\mathfrak{sl}_2}}

\DeclareMathOperator{\codim}{codim}
\DeclareMathOperator{\coker}{coker}

\begin{document}

\title[Discriminantal arrangements] 
{Resonant local systems on complements of discriminantal arrangements
and $\SL$ representations}

\author[D.~Cohen]{Daniel C.~Cohen$^\dag$}
\address{Department of Mathematics, Louisiana State University,
Baton Rouge, LA 70803}
\email{\href{mailto:cohen@math.lsu.edu}{cohen@math.lsu.edu}}
\urladdr{\href{http://www.math.lsu.edu/~cohen/}
{http://www.math.lsu.edu/\~{}cohen}}
\thanks{{$^\dag$}Partially supported by National Security Agency grant
MDA904-00-1-0038}

\author[A.~Varchenko]{Alexander N.~Varchenko$^\ddag$}
\address{Department of Mathematics, University of
North Carolina at Chapel Hill, Chapel Hill, NC 27599}
\email{\href{mailto:anv@email.unc.edu}{anv@email.unc.edu}}
\urladdr{\href{http://www.math.unc.edu/Faculty/av/}
{http://www.math.unc.edu/Faculty/av/}}
\thanks{{$^\ddag$}Partially supported by National Science 
Foundation grant DMS-9801582}

% 17B10 Representations, algebraic theory (weights)
%       (Lie algebras and Lie superalgebras)
% 32S22 Relations with arrangements of hyperplanes
%       (Several complex variables and analytic spaces; Singularities)
% 52C35 Arrangements of points, flats, hyperplanes
%       (Convex and discrete geometry; Discrete geometry)
% 55N25 Homology with local coefficients, equivariant cohomology
%       (Algebraic topology; Homology and cohomology theories)


\begin{abstract}
We calculate the skew-symmetric cohomology of the complement of a
discriminantal hyperplane arrangement with coefficients in local
systems arising in the context of the representation theory of the Lie
algebra $\SL$.  For a discriminantal arrangement in $\C^k$, the
skew-symmetric cohomology is nontrivial in dimension $k-1$ precisely
when the ``master function'' which defines the local system on the
complement has nonisolated critical points.  In symmetric coordinates,
the critical set is a union of lines.  Generically, the dimension of
this nontrivial skew-symmetric cohomology group is equal to the number
of critical lines.
\end{abstract}
\maketitle

\section{Introduction} \label{sec:intro}

Let $z=(z_1,\dots,z_n)$ be an $n$-tuple of distinct complex numbers,
$z_i \neq z_j$ for $i \neq j$, and let $m=(m_1,\dots,m_n)$ be an
$n$-tuple of nonnegative integers.  The ``master function''
\begin{equation} \label{eqn:master}
\Phi_{k,n}=\Phi_{k,n}(t;z,m) = \prod_{i=1}^k \prod_{j=1}^n 
(t_i - z_j)^{-m_j} \prod_{1\le p < q \le k}(t_p - t_q)^2
\end{equation}
corresponding to $z$ and $m$ arises in a number of contexts.  For
instance, hypergeometric integrals involving this function are used in
\cite{SV} to construct solutions of the $\SL$ KZ differential
equations.  Furthermore, the critical point equations of the master
function coincide with the Bethe ansatz equations for the $\SL$ Gaudin
model, see \cite{RV, V}.

Let $\A_{k,n}$ be the discriminantal arrangement in $\C^k$ consisting
of the hyperplanes $H^j_i=\{t_i=z_j\}$ and $H_{p,q}=\{t_p=t_q\}$
defined by the linear polynomials occuring in $\Phi_{k,n}$.  Denote
the complement of $\A_{k,n}$ by 
$X=X(\A_{k,n})=\C^k \setminus_{H \in \A_{k,n}} H$.  For 
$\kappa \in \C^*$, the function $\Phi_{k,n}^{1/\kappa}$ defines a
local system of coefficients $\LL$ on $X$, with monodromy
$\exp(2\pi\ii m_j/\kappa)$ about the hyperplane $H^j_i$, and monodromy
$\exp(-4\pi\ii/\kappa)$ about $H_{p,q}$.  These local systems are
often resonant, in the sense that the monodromy about intersections of
hyperplanes of $\A_{k,n}$ may be trivial.  For applications in
conformal field theory, the case where $\kappa$ is a positive
integer (greater than $2$) is of most interest.  In this instance, 
the aforementioned hypergeometric integrals are integrals of algebraic 
functions, providing motivation for the study of these local systems.

The symmetric group $\Sigma_k$ acts on $\C^k$ by permuting
coordinates.  This action preserves the complement $X$ of $\A_{k,n}$,
and also the local system $\LL$, as the monodromy about
$H^j_{\sigma(i)}$ (resp., $H_{\sigma(p),\sigma(q)}$) is the same as
that about $H^j_i$ (resp., $H_{p,q}$) for all $\sigma \in \Sigma_k$. 
Consequently, $\Sigma_k$ acts on the local system cohomology
$H^*(X;\LL)$.  Let $H^*_-(X;\LL)$ denote the subspace of all
skew-symmetric cohomology classes, those classes $\eta$ for which
$\sigma(\eta)=(-1)^{|\sigma|}\eta$ for all $\sigma \in \Sigma_k$.  
The purpose of this note is to determine the dimensions of these
skew-symmetric cohomology groups.

This determination is given in terms of the representation theory of
the Lie algebra $\SL$.  Let $e$, $f$, and $h$ be the standard
generators of $\SL$, satisfying $[e,f]=h$, $[h,e]=2e$, $[h,f]=-2f$. 
Let $L_a$ be the irreducible $\SL$-module with highest weight $a\in
\C$.  The module $L_a$ is generated by its singular vector $v_a$,
which satisfies $ev_a=0$, $hv_a=av_a$.  The vectors $v_a, fv_a,
f^2v_a\,,\dots $ form a basis for $L_a$.  If $a$ is a nonnegative
integer, then $\dim L_a=a+1$; otherwise $L_a$ is infinite-dimensional.
 
For nonnegative integers $m_1,\dots,m_n$ as above, the tensor product
$L^{\otimes m}=L_{m_1}\otimes\cdots \otimes L_{m_n}$ is a direct sum
of irreducible representations with highest weights $|m|-2k$, where
\[
|m| =m_1 + \cdots + m_n
\]
and $k$ is a nonnegative integer.  Let $w(m,k)$ denote the
multiplicity of $L_{|m|-2k}$ in $L^{\otimes m}$.  Then
\[
w(m,k)\geq 0\quad \text{if}\quad |m|-2k\geq 0 \qquad \text{and} 
\qquad w(m,k)=0\ \ {\rm if}\ \ |m|-2k<0.
\]

\begin{thm} \label{thm:dim}
Let $m_1,\dots,m_n$ be nonnegative integers, and let $\LL$ be the
local system on the complement $X$ of the discriminantal arrangement
$\A_{k,n}$ induced by $\left(\Phi_{k,n}(t;z,m)\right)^{1/\kappa}$ for
generic $\kappa$.
\begin{enumerate}
\item \label{item:dim1}
If $0\le |m|-k+1<k$, then $H^q_-(X;\LL)=0$ for $q< k-1$, 
$\dim H^{k-1}_-(X;\LL) = w(m,|m|-k+1)$, and 
$\dim H^k_-(X;\LL)=w(m,|m|-k+1)+\binom{n+k-2}{k}$.
\item  \label{item:dim2}
Otherwise, $H^q_-(X;\LL)=0$ for $q\neq k$ and 
$\dim H^k_-(X;\LL) = \binom{n+k-2}{k}$.
\end{enumerate}
\end{thm}

In \cite{ScV}, I.~Scherbak and the second author relate the critical
set of the master function $\Phi_{k,n}(t;z,m)$ of \eqref{eqn:master}
and the representation theory of the Lie algebra $\SL$ in the context
of Fuschian differential equations.  These results prompted us to
investigate the behavior of local systems on the complement of the
underlying discriminantal arrangement $\A_{k,n}$ induced by powers of
the master function, our aim being to determine if the aforementioned
relationship is reflected in the dimensions of the (skew-symmetric)
local system cohomology groups.  For the sake of comparison, we
briefly state results from \cite[Theorem 1]{ScV} concerning the
critical set of the function $\Phi_{k,n}(t;z,m)$.

\begin{thm} \label{thm:crit}
Let $m_1,\dots,m_n$ be nonnegative integers.
\begin{enumerate}
\item \label{item:crit1}
If $0\le |m|-k+1<k$, then for generic $z$, the function
$\Phi_{k,n}(t;z,m)$ has only nonisolated critical points.  In
symmetric coordinates $\lambda_1=\sum t_i, \ \lambda_2 = \sum t_i t_j,
\ \dots,\ \lambda_k=t_1\cdots t_k$, the critical set consists of
$w(m,|m|-k+1)$ lines. 
\item \label{item:crit2}
If $|m|-k+1 > k$, then for generic $z$, all critical points of
$\Phi_{k,n}(t;z,m)$ are nondegenerate and the critical set consists of
$w(m,k)$ $\Sigma_k$-orbits. 
\item  \label{item:crit3}
If $|m|-k+1=k$ or $|m|-k+1<0$, then for any $z$, the function
$\Phi_{k,n}(t;z,m)$ has no critical points in $X$.
\end{enumerate}
\end{thm}

Thus, for generic $z$, the skew-symmetric cohomology group
$H^{k-1}_-(X;\LL)$ is nontrivial if and only if the master function
$\Phi_{k,n}(t;z,m)$ has critical lines.  Moreover, the dimension of
$H^{k-1}_-(X;\LL)$ is equal to the number of critical lines.  Note
also that the behavior of the critical set of $\Phi_{k,n}(t;z,m)$
differs in cases~\ref{item:crit2} (isolated critical points) and
\ref{item:crit3} (no critical points) of Theorem \ref{thm:crit}, while
the skew-symmetric cohomology $H^{q}_-(X;\LL)$ vanishes in all
dimensions except $q=k$ in both cases.  It would be interesting to
determine if these cases are reflected in cohomological properties of
the local system $\LL$.

For an arbitrary arrangement $\A$ of $N$ hyperplanes in $\C^k$, the
set of complex rank one local systems on the complement $X(\A)$ may be
realized as the complex torus $\T\cong H^1(X(\A);\C^*)$.  
The correspondence is given by $\LL \leftrightarrow 
\tau=(\dots \ \tau_H\ \dots)\in \T$, where $\tau_H$ 
is the monodromy of $\LL=\LL(\tau)$ about the hyperplane $H$ of $\A$.  
Generically, the local system cohomology vanishes, except 
possibly in the top dimension, $H^q(X(\A);\LL)=0$ if
$q \neq k$.  Call these local systems nonresonant.
The local systems for which the cohomology does not vanish 
(for $q \neq k$) are called resonant, and correspond to 
elements of the cohomology jumping loci,
$V_p^q(\A) = \{\tau \in \T \mid \dim H^q(X(\A);\LL(\tau)) \ge p\}$,
which are known to be are unions of torsion-translated subtori of $\T$.  

The first cohomology jumping loci, $V^1_p(\A)$, have been the subject
of a great deal of recent attention, and are to some extent
understood.  There are, for instance, combinatorial algorithms for
determining the components passing through the identity in $\T$, see
\cite{CS, LY}.  Less is known about the higher jumping loci
$V^q_p(\A)$ for $1< q \le k$.  In the case where $\A=\A_{k,n}$ is a
discriminantal arrangement, Theorem \ref{thm:dim}.\ref{item:dim1}
provides new examples of resonant local systems, that is, nontrivial
elements of the varieties $V_p^q(\A_{k,n})$ for $q=k-1,k$, 
$1 \le p \le w(m,|m|-k+1)$, and arbitrary $k$.

This note is organized as follows.  Some results on the local system
cohomology of the complement of an arrangement, including a
strengthening of a particular case of 
Theorem \ref{thm:dim}.\ref{item:dim2}, are given in Section
\ref{sec:cohomology}.  See \cite{OT1,OT2} as general references in
this context.  In Section \ref{sec:skew}, we discuss the relationship
between the skew-symmetric local system cohomology of the complement
of a discriminantal arrangement and the representation theory of
$\SL$, and reformulate Theorem \ref{thm:dim} in terms of the latter in
Theorem \ref{thm:ker}.  After a number of preliminary results are
established in Section \ref{sec:prelim}, the proof of Theorem
\ref{thm:ker} is given in Section \ref{sec:proof}.

\section{Local system cohomology}
\label{sec:cohomology}

Choose coordinates $\b{t}=(t_1,\dots,t_k)$ on $\C^k$, and let $\A$ be
an arbitrary arrangement of $N$ hyperplanes in $\C^k$, with complement
$X=X(\A)=\C^k\setminus \bigcup_{H \in \A} H$.  Assume that $\A$
contains $k$ linearly independent hyperplanes.  For each hyperplane
$H$ of $\A$, let $f_H$ be a linear polynomial with
$H = \{\b{t} \in \C^k \mid f_H(\b{t})=0\}$, and let 
$\omega_H=d\log f_H$ denote the corresponding logarithmic one-form. 
Let $\la=(\dots \ \la_H\ \dots)\in\C^N$ be a weight vector, where
$\la_H\in\C$ corresponds to $H\in\A$.  Associated to $\la$, we have
\begin{enumerate}
\item[(1)] a flat connection on the trivial line bundle over $X$, with
connection form $\nabla=d+\omega_{\la}\wedge:\Omega^0 \to \Omega^1$,
where $d$ is the exterior differential operator with respect to the
coordintes $\b{t}$, $\omega_{\la} = \sum_{H\in\A} \la_H\,\omega_H$,
and $\Omega^q$ is the sheaf of germs of holomorphic differential forms
of degree $q$ on $X$;
\item[(2)] a rank one representation $\rho:\pi_1(X) \to \C^*$, given
by $\rho(\gamma_H)=\tau+H=\exp(-2\pi\ii\la_H)$, where $\gamma_H$ is any
meridian loop about the hyperplane $H$ of $\A$; and
\item[(3)] a rank one local system $\LL=\LL_{\la}$ on $X$ associated
to the flat connection $\nabla$ (resp., the representation $\rho$).
\end{enumerate}

Let $A(\A)$ be the Orlik-Solomon algebra of $\A$, generated by one
dimensional classes $a_H$ corresponding to the hyperplanes of $\A$. 
It is the quotient of the exterior algebra generated by these classes
by a homogeneous ideal, hence a finite dimensional graded
$\C$-algebra.  The weight vector $\la$ determines an element $a_\la
=\sum_{H\in\A} \la_H\, a_H\in A^1(\A)$.  Since $a_\la \wedge a_\la=0$,
multiplication by $a_\la$ defines a differential on $A(\A)$.  The
resulting complex $(A^\bul(\A),a_\la\wedge)$ may be identified with a
subcomplex of the twisted de Rham complex of $X$ with coefficients in
$\LL$ via $a_H \mapsto \omega_H$.

An edge of an arrangement $\A$ is a nonempty intersection of
hyperplanes in $\A$.  Associated to each flag $F$ of edges of $\A$ and
the weight vector $\la$, there is an element of the Orlik-Solomon
algebra and a corresponding logarithmic flag form $\varOmega_F$.  If
$F=(F_q \subset F_{q-1} \subset \cdots \subset F_2 \subset F_1)$ is a
flag of edges of $\A$ with $\codim F_j=j$ for each $j$, then
\begin{equation} \label{eqn:flag form}
\varOmega_F = \la_{F_1}\omega_{F_1} \wedge
\sum_{F_2 \subset H} \la_H \omega_H\wedge \cdots
\wedge \sum_{F_q \subset H} \la_H \omega_H.
\end{equation}

An edge is dense if the subarrangement of hyperplanes containing it is
irreducible: the hyperplanes cannot be partitioned into nonempty sets
so that, after a change of coordinates, hyperplanes in different sets
are in different coordinates.  Let $\Ai=\A\cup H_{\infty}$ be the
projective closure of $\A$, the union of $\A$ and the hyperplane at
infinity in $\CP^k$.  Set $\la_{H_{\infty}} = -\sum_{H \in \A}\la_H$,
and $\la_L=\sum_{L \subseteq H}\la_H$ for an edge $L$.

\begin{thm}[\cite{STV}]  \label{thm:STV}
Let $\LL$ be the rank one local system on the complement $X$ of $\A$
corresponding to the weight vector $\la$.
\begin{enumerate}
\item \label{item:comb}
If $\la_L \notin \Z_{>0}$ for every dense edge $L$ of $\Ai$, then
$H^*(X;\LL) \cong H^*(A^\bul(\A),a_\la\wedge)$.
\item \label{item:nonres}
If $\la_L \notin \Z_{\ge 0}$ for every dense edge $L$ of $\Ai$, then
$H^q(X;\LL) = 0$ for $q\neq k$ and $\dim H^k(X;\LL) = |\chi(X)|$,
where $\chi(X)$ is the Euler characteristic of $X$.
\end{enumerate}
\end{thm}

\begin{rem}\label{rem:equiv}
Call two weight vectors $\la$ and $\mu$ equivalent, and write 
$\la \equiv \mu$, if $\la - \mu$ is an integer vector.  Note that if
$\la\equiv\mu$ then $\exp(-2\pi\ii\la_j)= \exp(-2\pi\ii\mu_j)$ for each
$j$, so $\la$ and $\mu$ give rise to the same rank one local system
$\LL$ on $X$.  Consequently, if $\LL$ is the local system
corresponding to $\la$, and $\la$ is equivalent to a weight vector
which satisfies the conditions of 
Theorem \ref{thm:STV}.\ref{item:nonres} above, then $H^q(X;\LL)=0$ for
$q\neq k$.
\end{rem}

Now let $\A=\A_{k,n}$ be a discriminantal arrangement, and $\la$ the
weight vector corresponding to the master function $\Phi_{k,n}(t;z,m)$
of \eqref{eqn:master} and $\kappa\in \C^*$.  Explicitly, the weight of
the hyperplane $H^j_i$ is $\la^j_i=-m_j/\kappa$, and the weight of
$H_{i,j}$ is $\la_{i,j}=2/\kappa$.  Denote the corresponding local
system by $\LL_\kappa$ to indicate the dependence on $\kappa$.  The
following is an immediate consequence of 
Theorem \ref{thm:STV}.\ref{item:comb}.

\begin{prop}\label{prop:OS}
For generic $\kappa$, we have $H^*(X(\A_{k,n});\LL_\kappa)\cong
H^*(A^\bul(\A_{k,n}),a_\la\wedge)$.
\end{prop}

\begin{rem} \label{rem:res}
For an arbitrary arrangement $\A$ of $N$ hyperplanes in $\C^k$, the
resonance varieties of $\A$ are the cohomology jumping loci, 
$R_p^q(\A)= 
\{\la \in \C^N \mid \dim H^q(A^\bul(\A),a_\la\wedge) \ge p\}$, 
of the Orlik-Solomon complex, see \cite{Fa}.  
The variety $R_p^q(\A)$ coincides with the tangent cone of the 
cohomology jumping locus $V_p^q(\A)$ at the identity in $\T$, and 
is consequently a union of linear subspaces of 
$\C^N$, see for instance \cite{CO}.

Explicit combinatorial descriptions of the first resonance varieties
$R^1_p(\A)$ may be found in \cite{Fa,LY}.  Less is known about the
higher resonance varieties $R_p^q(\A)$ for $1 < q \le k$.  In the case
where $\A=\A_{k,n}$ is a discriminantal arrangement, and the weights
$\la$ correspond to the local system $\LL_\kappa$ for generic
$\kappa$, the local system cohomology is isomorphic to that of the
Orlik-Solomon complex by Proposition~\ref{prop:OS}.  Hence, Theorem
\ref{thm:dim}.\ref{item:dim1} provides nontrivial elements of the
varieties $R_p^q(\A_{k,n})$ for $q=k-1,k$, $1 \le p \le w(m,|m|-k+1)$,
and arbitrary $k$.
\end{rem}

The next result strengthens Theorem \ref{thm:dim}.\ref{item:dim2} 
in the case $|m|-k+1 \ge k$.

\begin{thm} \label{thm:nonres3}
If $|m|-k+1 \ge k$ and $\kappa$ is generic, then
$H^q(X(\A_{k,n});\LL_\kappa)=0$ for $q\neq k$ and 
$\dim H^k(X(\A_{k,n});\LL_\kappa)=(n+k-2)!/(n-2)!$
\end{thm}
\begin{proof}
Note that $|\chi(X(\A_{k,n}))|=(n+k-2)!/(n-2)!$.

If $k=1$, then $\A_{1,n}$ is an arrangement of $n$ points in $\C$.  In
this case, the condition $|m|-k+1 \ge k$ insures that $m_j \neq 0$ for
some $j$, $1\le j \le n$.  For any $\kappa$ for which $m_j/\kappa$ is
not an integer, the local system $\LL_\kappa$ on $X(\A_{1,n})$ is
nontrivial, and $\dim H^1(X(\A_{1,n});\LL_\kappa)=n-1$.

For $k\ge 2$, by Theorem \ref{thm:STV}.\ref{item:nonres} and Remark
\ref{rem:equiv}, it suffices to show that there is a weight vector
$\mu$ equivalent to $\la$ for which $\mu_L \notin \Z_{\ge 0}$ for
every dense edge $L$ of $\Ai$, where $\Ai$ is the projective closure
of the discriminantal arrangement $\A_{k,n}$.  We will show that there
are integers $a_1,\dots,a_n$, so that the weight vector $\mu$ given by
$\mu^j_i=a_j-\la^j_i=a_j-m_j/\kappa$ and
$\mu_{i,j}=\la_{i,j}=2/\kappa$ satisfies these conditions.

Denote the hyperplanes of $\Ai$ by $H^j_i$, $H_{i,j}$, and $H_\infty$. 
The dense edges of $\Ai$ may be described as follows.  For $I\subseteq
[k]$ and $j \in [n]$, let $L_I= \bigcap_{p,q\in I} H_{p,q}$ and
$L^j_I=\bigcap_{i\in I} H^j_i$.  If $|I|=1$, set $L_I =\C^k$.  For
$\emptyset \neq J \subseteq [n]$, set $L^J_I=\bigcap_{j\in J} L^j_I$. 
One can check that the dense edges of $\Ai$, and their weights with
respect to $\mu$, are then given by
\[
\begin{matrix}
&\text{dense edge $L$} \hfill &\text{weight $\mu_L$} \hfill \\
\\
\text{(a)} & L_I,\ I \subseteq [k],\ |I|=l,\ 2 \le l \le k \hfill &
l(l-1)/\kappa \hfill \\
\\
\text{(b)} & L^j_I \cap L_I,\ j \in [n],\ I \subseteq [k],\ |I|=l
\hfill & la_j +l(l-m_j-1)/\kappa \hfill \\
\\
\text{(c)} & H_\infty \hfill& k(|m|-k+1)/\kappa - k|a| \hfill \\
\\
\text{(d)} & H_\infty \cap L_i^{[n]},\ i \in [k] \hfill &
(k-1)(|m|-k)/\kappa - (k-1)|a| \hfill \\
\\
\text{(e)} & H_\infty \cap L^J_I \cap L_I,\ I \subseteq [k],\ |I|=l,\
2 \le l < k \quad & (k-l)(|m|-k+1)/\kappa + l|m^J|/\kappa - k|a| +
l|a_J| \\
& \phantom{H_\infty \cap L^J_I \cap L_I,\ }J \subseteq [n],\ |J| \ge 2
\hfill
\end{matrix}
\]
In (e) above, we use the notation $|m^J| = |m| - \sum_{j \in J} m_j$
and $|a_J|=\sum_{j\in J}a_j$.

Now assume that $|m|-k+1 \ge k$ and that $\kappa$ is generic.  These
conditions imply that the weights in (a), (c), (d), and (e) above are
not integers, for any choice of $a=(a_1,\dots,a_n)$.  Choosing
$a_j=-1$ whenever $1\le m_j \le k-1$ insures that the weights in (b)
are not in $\Z_{\ge 0}$.  The result follows.
\end{proof}

\section{Skew-symmetric cohomology and $\SL$ representations}
\label{sec:skew}

We now turn to the skew-symmetric cohomology groups
$H^q_-(X(\A_{k,n});\LL_\kappa)$, and their relation to representations
of $\SL$.

For $a\in\C$, let $M_a$ be the corresponding Verma module, the
infinite dimensional $\SL$-module generated by the vector $v_a$, where
$e v_a=0$ and $h v_a=a v_a$.  The vectors $\{f^kv_a\mid k \ge 0\}$
form a basis for $M_a$, and the $\SL$ action is given by
\[
e:f^kv_a \mapsto k(m-k+1)f^{k-1}v_a, \qquad
f:f^kv_a \mapsto f^{k+1}v_a, \qquad
h:f^kv_a \mapsto (m-2k)f^kv_a.
\]
Recall \cite{Kac} that the Shapovalov form is the symmetric bilinear
form $S_a$ on $M_a$ defined by $S_a(v_a,v_a)=1$, $S_a(hx,y)=S_a(x,hy)$
and $S_a(fx,y)=S_a(x,ey)$ for all $x,y \in M_a$, and
\[
S(f^k v_a,f^l v_a)=\begin{cases}
k!\, a(a-1) \cdots (a-k+1) &\text{if $k=l$,}\\
0 &\text{if $k \neq l$.}
\end{cases}
\]

Given $m=(m_1,\dots,m_n)\in\Z^n$, let 
$M^{\otimes m}=M_{m_1} \otimes \cdots \otimes M_{m_n}$ denote the
corresponding tensor product of Verma modules.  A basis for
$M^{\otimes m}$ is given by
\[
F_Jv:= f^{j_1}v_{m_1}\otimes\cdots \otimes f^{j_n}v_{m_n}, \quad
\text{where $J=(j_1,\dots,j_n)$ with $j_i \ge 0$ for each $i$.}
\]
The action of $\SL$ on $M^{\otimes m}$ is given by
\begin{equation} \label{eqn:action}
e:F_Jv  \mapsto \sum_{i=1}^n j_i (m_i - j_i + 1) F_{J-1_i}v, \quad
f:F_Jv  \mapsto \sum_{i=1}^nF_{J+1_i}v, \quad
h:F_Jv  \mapsto (|m|-2|J|)F_Jv,
\end{equation}
where $J \pm 1_i = (j_1,\dots,j_i\pm 1,\dots,j_n)$.  
Let $S$ denote the Shapovalov form 
$S_{m_1} \otimes \cdots \otimes S_{m_n}$ 
on $M^{\otimes m}$.

For an $\SL$-module $V$ and $\lambda\in\C$, let 
$V[\lambda]=\{x \in V \mid hx=\lambda x\}$ be the weight subspace of
weight $\la$.  For the tensor product $M^{\otimes m}$ and an integer
$k$ such that $|m|-2k \ge 0$, the weight subspace 
$M^{\otimes m}[|m|-2k]$ has basis
\[
F_Jv:= f^{j_1}v_{m_1}\otimes \cdots \otimes f^{j_n}v_{m_n},
\]
where $J$ runs through all multiindices such that 
$|J|=j_1 +\dots +j_n = k$ with nonnegative $j_i$.  
The dual space $(M^{\otimes m}[\|m|-2k])^*$ has the dual basis, 
denoted by 
\[
(F_Jv)^*:= (f^{j_1}v_{m_1}\otimes \cdots \otimes f^{j_n}v_{m_n})^*.
\]
Note that $\dim (M^{\otimes m}[|m|-2k])^* = \binom{n+k-1}{k}$.

Let $(M^{\otimes m})^*=\bigoplus_k (M^{\otimes m}[|m|-2k])^*$ denote
the restricted dual of the $\SL$-module $M^{\otimes m}$, with basis
$(F_Jv)^*$ as above, for all relevant $J$.  The contragredient action
of $\SL$ on $(M^{\otimes m})^*$ is given by
\begin{equation} \label{eqn:contra}
\begin{aligned}
f: (F_Jv)^*&\mapsto \sum_{i=1}^n (j_i+1)(m_i-j_i)(F_{J+1_i}v)^*,\\
e:(F_Jv)^*&\mapsto \sum_{j_i \neq 0}(F_{J-1_i}v)^*,\\
h:(F_Jv)^*&\mapsto (|m|-2|J|)(F_Jv)^*.
\end{aligned}
\end{equation}
The Shapovalov form gives rise to a homomorphism of $\SL$-modules
$S:M^{\otimes m} \to (M^{\otimes m})^*$ defined by 
\begin{equation} \label{eqn:Shapovalov}
S(F_Jv) = c_J(F_Jv)^*,\quad \text{where}
\quad c_J=\prod_{i=1}^n j_i!\, m_i (m_i - 1) \cdots 
(m_i - j_i + 1).
\end{equation}

By Proposition \ref{prop:OS}, the skew-symmetric cohomology
$H^*_-(X(\A_{k,n});\LL_\kappa)$ is isomorphic to the skew-symmetric
part of the cohomology of the Orlik-Solomon complex
$(A^\bul(\A_{k,n}),a_\la\wedge)$.  By results of \cite{SV}, this in
turn is isomorphic to the cohomology of the complex
\begin{equation} \label{eqn:complex}
0 \to (M^{\otimes m}[ |m| -2k + 2])^* \xrightarrow{\ \ f\ \ }
( M^{\otimes m}[ |m| -2k])^* \to 0. 
\end{equation}
Recall that the Orlik-Solomon complex may be realized as a subcomplex
of the twisted de Rham complex of $X(\A_{k,n})$ with coefficients in
$\LL_\kappa$ by identifying generators with the corresponding
logarithmic forms.  Recall also that the hyperplanes of $\A_{k,n}$
include $H_i^j=\{t_i-z_j=0\}$.  

Associated to each monomial 
$(F_Jv)^*=(f^{j_1}v_{m_1}\otimes \cdots \otimes f^{j_n}v_{m_n})^*$ in
$(M^{\otimes m})^*$, there is a skew-symmetric logarithmic form
$\omega_J$ defined as follows.  Given $(F_Jv)^*$, 
let $\ell_i(J)=j_1+\cdots +j_i$ for $1\le i \le n$.  Note that 
$\ell_n(J)=|J|$, and set $\ell_0(J)=0$.
Define
$\eta_J =  \alpha_J\, \eta_{J,1}\wedge \eta_{J,2} \wedge \cdots
\wedge \eta_{J,n}$, 
where
\[
\alpha_J = \frac{1}{j_1!\, j_2! \cdots j_n!} \quad\text{and}\quad
\eta_{J,i} =  \frac{d (t_{\ell_{i-1}(J)+1} -z_i)}{t_{\ell_{i-1}(J)+1} -z_i} 
\wedge \cdots 
\wedge
\frac{d (t_{\ell_i(J)-1} -z_i)}{t_{\ell_i(J)-1} -z_i}
\wedge
\frac{d (t_{\ell_i(J)} -z_i)}{t_{\ell_i(J)} -z_i},
\]
and let $\omega_J$ be the skew-symmetrization of $\eta_J$. 
Identifying $\omega_J$ with the corresponding element of
$A(\A_{k,n})$, this defines a map from $(M^{\otimes m}[|m|-2j])^*$ to
the Orlik-Solomon complex,\begin{equation} \label{eqn:skew
form}(F_Jv)^* \mapsto \omega_J,\end{equation} for each $j$, which
induces an isomorphism between the cohomology of the complex
\eqref{eqn:complex} and the skew-symmetric cohomology of
$(A^\bul(\A_{k,n}),a_\la\wedge)$.  

The complex \eqref{eqn:complex} is located in dimensions $k-1$ and
$k$.  Consequently, $H^{k-1}_-(X(\A_{k,n});\LL_\kappa)\cong \ker f$,
$H^{k}_-(X(\A_{k,n});\LL_\kappa)\cong \coker f$, and the skew-symmetric
cohomology groups vanish in other dimensions,
$H^q_-(X(\A_{k,n});\LL_\kappa)=0$ for $q \neq k-1,k$.  The
differential $f$ of the complex \eqref{eqn:complex} is given by the
action of $f\in\SL$ on $(M^{\otimes m})^*$ recorded in
\eqref{eqn:contra} above.  Recall from the Introduction that
$L^{\otimes m} = L_{m_1} \otimes \cdots \otimes L_{m_n}$ is the tensor
product of the irreducible $\SL$-modules $L_{m_i}$, $1\le i \le n$,
and that $w(m,j)$ denotes the multiplicity of $L_{|m|-2j}$ in
$L^{\otimes m}$.

\begin{thm} \label{thm:ker}
Let $m_1,\dots,m_n \in \Z_{\ge 0}$, and 
$f:(M^{\otimes m}[ |m| -2k + 2])^* \to( M^{\otimes m}[ |m| -2k])^*$.
\begin{enumerate}
\item \label{item:ker1}
If $0\le |m|-k+1<k$, then the dimension of the kernel of $f$ is 
$\dim \ker f = w(m,|m|-k+1)$, and 
$\dim \coker f=w(m,|m|-k+1)+\binom{n+k-2}{k}$.
\item  \label{item:ker2}
Otherwise, $\ker f=0$ and 
$\dim \coker f= \binom{n+k-2}{k}$.
\end{enumerate}
\end{thm}

\begin{rem} \label{rem:ker&dim}
Since the skew-symmetric cohomology $H^*_-(X(\A_{k,n});\LL_\kappa)$ 
is isomorphic to the cohomology of the complex \eqref{eqn:complex},
Theorem \ref{thm:dim} from the Introduction is a consequence of
Theorem \ref{thm:ker}.  Consequently, the next two sections are
devoted to establishing this latter result.
\end{rem}

\begin{rem} \label{rem: big k}
If $m_i \ge k$ for some $i$, $1\le i \le n$, then $\ker f = 0$ in
either case of Theorem \ref{thm:ker} above, 
see Theorem \ref{thm:k2}.\ref{item:k21}. 
Accordingly, in the case $0 \le |m| - k + 1 < k$ and $m_i \ge k$ for
some $i$, we have $w(m,|m|-k+1)=0$, see Lemma \ref{lem:consistent}.
\end{rem}

For the tensor product $L^{\otimes m}=L_{m_1}\otimes\cdots\otimes
L_{m_n}$, a basis for the weight subspace $L^{\otimes m}[|m|-2k]$ is
given by monomials$F_Jv =
f^{j_1}v_1\otimes\cdots \otimes f^{j_n}v_n$, where $j_1,\dots,j_n$ are
integers satisfying $j_1+\cdots+j_n=k$ and $0 \le j_i \le m_i$ for
each $i$.  The map \eqref{eqn:Shapovalov}
induced by the Shapovalov form gives rise to an injective map of
complexes
\[
\begin{CD}
0 @>>> L^{\otimes m}[|m|-2k+2] @>f>>  L^{\otimes m}[|m|-2k]  @>>> 0\\
@. @VVSV @VVSV  \\
0 @>>> (M^{\otimes m}[|m|-2k+2])^* @>f>>  
(M^{\otimes m}[|m|-2k])^* @>>> 0
\end{CD}
\]
where the action of $f$ is given by \eqref{eqn:action}  
on $L^{\otimes m}[|m|-2k+2]$, and by 
\eqref{eqn:contra} on $(M^{\otimes m}[|m|-2k+2])^*$.

As noted in the Introduction, the tensor product $L^{\otimes m}$ is a
direct sum of irreducible $\SL$-modules of the form $L_{|m|-2j}$.  For
each such summand $L_a$, the vector $f^av$ of lowest weight is in the
kernel of $f$.  This observation yields the following.

\begin{prop} \label{prop:flag}
The dimension
of the kernel the map
$f: L^{\otimes m}[ |m| -2k + 2] \to L ^{\otimes m}[ |m| -2k]$
is equal to $w(m, |m| + 1 - k)$ if $0\leq |m| -k + 1 < k$.
\end{prop}

\begin{rem} \label{rem:flag basis}
Via the embedding 
$S:L^{\otimes m}[ |m| -2k + 2] \to (M^{\otimes m}[|m|-2k+2])^*$
induced by the Shapovalov form, this result yields a subspace of
$\ker\bigl(f:(M^{\otimes m}[|m|-2k+2])^* \to 
(M^{\otimes m}[|m|-2k])^*\bigr)$ of dimension $w(m,|m|-k+1)$, with
basis corresponding to lowest weight vectors in the case 
$0\le |m|-k+1 <k$.  These, in fact, form a basis for
$\ker\bigl(f:(M^{\otimes m}[|m|-2k+2])^* \to 
(M^{\otimes m}[|m|-2k])^*\bigr)$, as asserted in Theorem
\ref{thm:ker}.\ref{item:ker1}, and shown in Section \ref{sec:proof}.

Identifying elements of $(M^{\otimes m})^*$ with skew-symmetric
logarithmic forms using \eqref{eqn:skew form}, the image of the
embedding 
$S:L^{\otimes m}[ |m| -2k + 2] \to (M^{\otimes m}[|m|-2k+2])^*$ is
realized by skew-symmetric flag forms, see \eqref{eqn:flag form}.  
Thus, upon establishing Theorem \ref{thm:ker}.\ref{item:ker1}, the
above considerations yield a basis of the skew-symmetric cohomology
group $H^{k-1}_-(X(\A_{k,n});\LL_\kappa)$ consisting of flag forms in
the case $0\le |m|-k+1 <k$.
\end{rem}

\section{Preliminary results} \label{sec:prelim}

The purpose of this section is to establish a number of results which
will be of use in the proof of Theorem \ref{thm:ker}.  First, we
describe the matrix of the endomorphism
\[
f:(M^{\otimes m}[ |m| -2k + 2])^* \to (M^{\otimes m}[ |m| -2k])^*.
\]
Order the bases of $(M^{\otimes m}[ |m| -2k + 2])^*$ and 
$(M^{\otimes m}[|m|-2k])^* $ lexicographically:  
$(F_Jv)^* > (F_Iv)^*$ if the left-most entry in $J-I$ is positive,
where $J=(j_1,\dots,j_n)$ and $I=(i_1,\dots,i_n)$.  Let $A_{k}(m)$
denote the matrix of $f:(M^{\otimes m}[ |m| -2k + 2])^* \to
(M^{\otimes m}[ |m| -2k])^*$ with respect to these ordered bases,
and acting on row vectors.  If $m=(m_1,m_2,\dots,m_n)$, let
$m^1=(m_2,\dots,m_n)$.

\begin{prop} \label{prop:matrix}
The matrix of $f:(M^{\otimes m}[ |m| -2k + 2])^* \to 
(M^{\otimes m}[|m| -2k])^*$ is given by
\[
A_{k}(m)=
\begin{pmatrix}
D_{1,k}(m) & A_{1}(m^1) & \b{0} & \cdots & \cdots & \b{0} \\
\b{0} & D_{2,k}(m) & A_{2}(m^1)  &&& \vdots \\
\vdots & & \ddots & \ddots & & \vdots \\
\vdots & &  & D_{k-1,k}(m) & A_{k-1}(m^1) & \b{0} \\
\b{0} & \cdots &\cdots & \b{0} & D_{k,k}(m) & A_{k}(m^1)
\end{pmatrix},
\]
where $A_{q}(m^1)$ is the matrix of $f:(M^{\otimes m^1}[|m^1|-2q+2])^*
\to (M^{\otimes m^1}[|m^1|-2q])^*$, $D_{q,k}(m)$ is the diagonal
matrix $(k-q+1)(m_1-k+q)\b{I}$, and $\b{0}$ and $\b{I}$ denote the
zero matrix and identity matrix of the appropriate sizes.
\end{prop}
\begin{proof}
Suppose the ordered basis of $(M^{\otimes m}[ |m| -2k + 2])^*$
corresponds to the $n$-tuples
\[
J_{1,1},\dots, J_{1,p_1}, J_{2,1},\dots, J_{2,p_2},\dots\dots, 
J_{n,1},\dots, J_{n,p_n},
\]
where $J_{i,1},\dots,J_{i,p_i}$ correspond to those basis elements for
which $j_1=\dots=j_{p-1} =0$ and $j_p\neq 0$.  Then the ordered basis
for $(M^{\otimes m}[ |m| -2k ])^*$ corresponds to the $n$-tuples
\[
\begin{aligned}
J_{1,1}+1_1,\dots, J_{1,p_1}+1_1, J_{2,1}+1_1,\dots, J_{2,p_2}+1_1, 
\dots\dots,\ 
&J_{n,1}+1_1,\dots, J_{n,p_n}+1_1,\\
J_{2,1}+1_2,\dots, J_{2,p_2}+1_2,\dots\dots,\ 
&J_{n,1}+1_2,\dots, J_{n,p_n}+1_2, \\
\ddots \quad & \phantom{J_{n,1}+1_n,\dots,J_{n,p_n}} \vdots \\
&J_{n,1}+1_n,\dots, J_{n,p_n}+1_n. 
\end{aligned}
\]
So if row $i$ of $A_{k}(m)$ corresponds to the basis element
$(F_Jv)^*$ of $(M^{\otimes m}[ |m| -2k + 2])^*$, then column $i$
corresponds to the basis element $(F_{J+1_1}v)^*$ of 
$(M^{\otimes m}[|m|-2k])^*$.  Hence, the diagonal entries are
$\left(A_{k}(m)\right)_{i,i}=(j_1+1)(m_1-j_1)$.  Furthermore, since
$J+1_1 > J+1_2 > \cdots > J+1_n$ in the lexicographic ordering, the
entries below the diagonal are $\left(A_{k,n}(m)\right)_{i,j}=0$ for
$i>j$.

Since
\[
\begin{aligned}
f((f^{j_1}v_{m_1}\otimes 
\cdots \otimes f^{j_n}v_{m_n})^*)=&(j_1+1)(m_1-j_1)(f^{j_1+1}v_{m_1}
\otimes 
f^{j_2}v_{m_2}\otimes \cdots \otimes f^{j_n}v_{m_n})^*\\
&\qquad+(f^{j_1}v_{m_1})^* \otimes f((f^{j_2}v_{m_2}\otimes 
\cdots \otimes f^{j_n}v_{m_n})^*),
\end{aligned}
\]
the fact that the (nonzero) entries above the diagonal are as asserted
also follows from the above considerations.
\end{proof}

\begin{exm} \label{exm:n=2}
In the case $n=2$, the matrix $A_k(m_1,m_2)$ 
has two nonzero entries in each row, and 
is given by
\[
\begin{pmatrix}
k(m_1-k+1) & m_2 \\
&& \ddots \hfill \\
& & (k-i+1)(m_1-k+i) & i(m_2-i+1) \\
&&&  \hfill \ddots \\
&&&& m_1 & k(m_2-k+1)
\end{pmatrix}.
\]
Using this, an exercise in linear algebra reveals that 
$\dim \ker A_2(m_1,m_2)=1$ if $0\le m_1,m_2<k$ and 
$0 \le m_1+m_2-k+1 < k$, and that $\ker A_k(m_1,m_2)=0$ otherwise.  It
is readily checked that this is the content of Theorem \ref{thm:ker}
in the case $n=2$.
\end{exm}

For any $n$, Proposition \ref{prop:matrix} has the following
consequence, which may be established using elementary row and column
operations.

\begin{cor} \label{cor:rowcolequiv}
Assume that $m_1=k-p$ for some $p$, $1\le p \le k$.  Then 
the matrix $A_{k}(m)$ of 
$f:(M^{\otimes m}[|m|-2k+2])^* \to (M^{\otimes m}[|m|-2k])^*$ 
is equivalent, via elementary row and column operations, to
\[
\begin{pmatrix}
\b{I} & A_{1}(m^1) \\
& \ddots & \ddots  \\
&& \b{I} & A_{p-1}(m^1)\\
&&&\b{0} & A_{p}(m^1)\\
&&&& \b{I}& A_{p+1}(m^1)\\
&&&&&\ddots&\ddots\\
&&&&&& \b{I} & A_{k}(m^1)
\end{pmatrix}.
\]
\end{cor}

Next, we record some properties of the multiplicity $w(m,j)$ of the
irreducible representation $L_{|m|-2j}$ in the tensor product
$L^{\otimes m}$.  Let $m=(m_1,m_2,\dots,m_n)$, and assume without loss
of generality that $0 \le m_1 \le m_2 \le \dots \le m_n$.  Recall that
$m^1=(m_2,\dots,m_n)$.  Let $r,r^1 \in \Z$ be maximal so that $|m|-2r
\ge 0$ and $|m^1|-2r^1 \ge 0$.  For $n \ge 3$, one can show that $r^1
\ge m_1$.

\begin{lem} \label{lem:recursion}
For $n \ge 3$, the multiplicity $w(m,j)$ of $L_{|m|-2j}$ in
$L^{\otimes m}$ satisfies
\[
w(m,j)=
\begin{cases}
w(m^1,0)+\dots+w(m^1,j-1)+w(m^1,j) 
& \text{if $0 \le j < m_1$,} \\
w(m^1,j-m_1)+\dots+w(m^1,j-1)+w(m^1,j) 
& \text{if $m_1 \le j \le r^1$,} \\
w(m^1,j-m_1)+\dots+w(m^1,|m^1|-j-1)+w(m^1,|m^1|-j) 
& \text{if $r^1 < j \le r$.}
\end{cases}
\]
\end{lem}
\begin{proof}
This is an elementary, albeit delicate, exercise using the fact that
if $a \le b$ are nonnegative integers, then $L_a \otimes L_b =
\bigoplus_{i=0}^a L_{a+b-2i}= L_{a+b} \oplus L_{a+b-2} \oplus \cdots
\oplus L_{b-a}$.
\end{proof}

\begin{lem} \label{lem:consistent}  
If $0 \le |m|-k+1 < k$ and there exists an $i$ for which $m_i \ge k$,
then $w(m,|m|-k+1)=0$.
\end{lem}
Note that if $|m|-k+1<k$, there can be at most one $i$ for which
$m_i\ge k$.
\begin{proof}
Let $m=(m_1,\dots,m_n)$.  The proof is by induction on $n$, with the
case $n=1$ trivial.  The case $n=2$ is also known to hold, as noted in
\cite{ScV}.

For general $n \ge 3$, assume that $m_1 \le m_2 \le \dots \le m_{n-1}
< k \le m_n$.  Take $j=|m|-k+1$, $0\le j < k$ in each of the three
cases in Lemma \ref{lem:recursion} above.  Write $m_1=k-p$, so
$|m^1|=|m|-m_1=|m|-k+p$.

If $0 \le j < m_1$, then
\[
|m^1|-p+1 = j,\ |m^1|-(p+1)+1 = j-1,\ \dots,\ |m^1|-(p+j)+1=0.
\]
Since $|m^1|-(p+j)+1=0$ and $j < m_1 \le r^1$, we have $|m^1|=p+j-1
\ge 2r^1 \ge 2m_1 \ge 2j$.  So $p-1 \ge j$ and $p>j$.  It follows that
$0\le |m^1|-(p+i)+1 < p+i$ for $i=0,1,\dots,j$.  So by induction,
$w(m^1,|m^1|-(p+i)+1)=w(m^1,j-i) =0$ for $i=0,1,\dots,j$.  Thus
$w(m,j)= w(m^1,j)+w(m^1,j-1)+\dots+w(m^1,1)+w(m^1,0)=0$ if 
$0 \le j <m_1$.

If $m_1 \le j \le r^1$, then
\[
|m^1|-p+1 = j,\ |m^1|-(p+1)+1 = j-1,\ \dots,\ |m^1|-(k-1)+1=j-m_1+1,\
|m^1|-k+1=j-m_1.
\]
are all nonnegative.  Since $|m^1|-k+1=j-m_1$ and $j \le r^1$, we
have $|m^1|=j-m_1+k-1 \ge 2r^1 \ge 2j$.  So $p-1 \ge j$ and $p > j$. 
It follows that $0\le |m^1|-(p+i)+1 < p+i$ for $i=0,1,\dots,k-p$.  So
by induction, $w(m^1,|m^1|-(p+i)+1)=w(m^1,j-i) =0$ for
$i=0,1,\dots,k-p$.  Thus $w(m,j)=
w(m^1,j)+w(m^1,j-1)+\dots+w(m^1,j-m_1+1)+w(m^1,j-m_1)=0$ if 
$m_1 \le j \le r^1$.

If $r^1 < j \le r$, then
\[
|m^1|-k+1 = j-m_1,\ |m^1|-(k-1)+1 = j-m_1+1,\ \dots,\
|m^1|-(k-|m|+2j)+1=|m^1|-j.
\]
Since $|m|-k+1=j$, we have $k+j-|m|=1$.  The condition $r^1 < j$
implies that $2j > |m^1|$.  It follows that $k-|m|+2 j-|m^1|+j>0$,
that is, $k-|m|+2j>|m^1|-j$.  From this, we obtain $0\le
|m^1|-(k-i)+1=j-m_1+i <k-i$ for $i=0,1,\dots,|m|-2j$.  So by
induction, $w(m^1,|m^1|-(k-i)+1)=w(m^1,j-m_1+i) =0$ for
$i=0,1,\dots,|m|-2j$.  Thus $w(m,j)=
w(m^1,j-m_1)+\dots+w(m^1,|m^1|-j-1)+w(m^1,|m^1|-j) =0$ 
if $r^1 < j \le r$.
\end{proof}

\section{Proof of Theorem \ref{thm:ker}} \label{sec:proof}

With the results of the previous sections at hand, 
we prove Theorem \ref{thm:ker}.  

\begin{thm} \label{thm:k2}
Let $m_1,\dots,m_n \in \Z_{\ge 0}$, and 
$f:(M^{\otimes m}[ |m| -2k + 2])^* \to( M^{\otimes m}[ |m| -2k])^*$.  
If 
\begin{enumerate}
\item \label{item:k21}
$m_i \ge k$ for some $i$, $1\le i \le n$, or 
\item \label{item:k22}
$|m|-k+1<0$, or 
\item  \label{item:k23}
$|m|-k+1\ge k$,
\end{enumerate}
then $\ker f=0$ and $\dim \coker f = \binom{n+k-2}{k}$.
\end{thm}
Note that Theorem \ref{thm:ker}.\ref{item:ker2} follows from this
result.
\begin{proof}
We will show that $\ker f=0$ in each of the three cases above
separately.  Proposition \ref{prop:matrix} facilitates elementary
proofs of cases \ref{item:k21} and \ref{item:k22}.

In case \ref{item:k21}, where $m_i \ge k$ for some $i$, 
we may assume without loss of generality that 
$m_1 \ge k$.  In this instance, each of the diagonal
matrices $D_{q,k}(m)=(k-q+1)(m_1-k+q)\b{I}$ occuring 
in the matrix $A_k(m)$ of $f$ is invertible. 
It follows immediately that $\ker f=0$.

In case \ref{item:k22}, where $|m|-k+1<0$, the proof is by double
induction on $k$ and $n$.  The result holds trivially for $n=1$ since
$A_{k}(m)=k(m_1-k+1)\neq 0$ if $|m|=m_1<k-1$.  For $m_i\ge 0$, the
condition $|m|-k+1 < 0$ is vaccuous if $k=1$, but the assertion holds
for $k=2$ and any $n$, as is readily checked.

In general, write $m_1=k-p$ for some $p$, $1\le p \le k$.  Then the
matrix $A_{k}(m)$ is of the form
\begin{equation} \label{eqn:matrix form}
\begin{pmatrix}
D_{1,k} & A_{1}(m^1) \\
& \ddots & \ddots  \\
&& D_{p-1,k} & A_{p-1}(m^1)\\
&&&\b{0} & A_{p}(m^1)\\
&&&&D_{p+1,k} & A_{p+1}(m^1)\\
&&&&&\ddots&\ddots\\
&&&&&& D_{k,k} & A_{k}(m^1)
\end{pmatrix}
\end{equation}
where $D_{q,k}=D_{q,k}(m)=(k-q+1)(m_1-k+q)\b{I}$ is nonzero if $q\neq
p$.  Since $m_1=k-p$ and $|m|-k+1<0$, we have $|m^1|-q+1<0$ for $p \le
q \le k$.  So $\ker A_{q}(m^1)=0$ for each such $q$ by induction.  The
result is obtained from these observations as follows.

Suppose $v=\begin{pmatrix}v_1 & \cdots & v_{p-1} & v_p & \cdots &
v_k\end{pmatrix}$ is in $\ker A_{k}(m)$.  Then from the invertibility
of $D_{q,k}(m)$ for $q < p$, we successively get $v_1=0,\ v_2=0,\
\dots,\ v_{p-1}=0$.  Similarly, from the invertibility of $D_{q,k}(m)$
for $q > p$ and the fact that $\ker A_{q}(m^1)=0$ for $p\le q\le k$,
we successively get $v_k=0,\ v_{k-1}=0,\ \dots,\ v_{p}=0$.

In case \ref{item:k23}, where $|m|-k+1\ge k$, recall from 
Section \ref{sec:skew} that the kernel of $f$ is isomorphic to
$H^{k-1}_-(X(\A_{k,n});\LL_\kappa)$, the skew-symmetric part of the
$(k-1)$-st cohomology of the complement of the discriminantal
arrangement $\A_{k,n}$ with coefficients in the local system induced
by $\left(\Phi_{k,n}(t;z,m)\right)^{1/\kappa}$ for generic $\kappa$. 
If $|m|-k+1 \ge k$, the local system cohomology
$H^{k-1}(X(\A_{k,n});\LL_\kappa)$ vanishes by 
Theorem \ref{thm:nonres3}.  Hence, in this instance, 
$\ker f \cong H^{k-1}_-(X(\A_{k,n});\LL_\kappa) \subseteq
H^{k-1}(X(\A_{k,n});\LL_\kappa)=0$.
\end{proof}

It remains to prove assertion \ref{item:ker1} of 
Theorem \ref{thm:ker}.

\begin{proof}[Proof of Theorem \ref{thm:ker}.\ref{item:ker1}]
Since $\dim (M^{\otimes m}[|m|-2j])^* = \binom{n+j-1}{j}$, 
we must show that
\[
\dim \ker\bigl(f:(M^{\otimes m}[|m|-2k+2])^* \to 
(M^{\otimes m}[|m|-2k])^*\bigr)
=\dim \ker A_k(m) = w(m,|m|-k+1)
\]
if $0 \le |m|-k+1 <k$.  By Proposition \ref{prop:flag} and Remark
\ref{rem:flag basis}, $\ker f = \ker A_k(m)$ contains a subspace of
dimension $w(m,|m|-k+1)$ generated by flag forms.  Thus, $\dim \ker
A_k(m) \ge w(m,|m|-k+1)$, and to establish the result, it suffices to
show that $\dim \ker A_k(m) \le w(m,|m|-k+1)$.

The proof is by double induction on $k$ and $n$.  The result holds for
$n=1,2$ and any $k$, and also for $k=2$ and any $n$ by direct
calculation, see Proposition \ref{prop:matrix} and Example
\ref{exm:n=2}.  (The condition $0\le |m|-k+1 < k$ is vacuous for
$k=1$.)  So assume that $k \ge 3$ and $n \ge 3$.

If $m_1 \ge k$, then $\ker f = 0$ by 
Theorem \ref{thm:k2}.\ref{item:k21} and $w(m,|m|-k+1)=0$ by Lemma
\ref{lem:consistent}.  So assume that $0\le m_1 < k$, and write
$p=k-m_1$.  Then the matrix $A_{k}(m)$ of 
$f:(M^{\otimes m}[|m|-2k+2])^* \to (M^{\otimes m}[|m|-2k])^*$ is as
given in the proof of Theorem \ref{thm:k2}.\ref{item:k22}, see
\eqref{eqn:matrix form}.  Since the diagonal matrix 
$D_{q,k}=D_{q,k}(m)=(k-q+1)(k-q-m_1)\b{I}$
is invertible for $q\neq p$, any nontrivial
element of $\ker A_k(m)$ is necessarily of the form
\begin{equation} \label{eqn:kernel elements}
v=\begin{pmatrix}
0 & \cdots & 0 & v_p & v_{p+1} & \cdots & v_q & 0 & \cdots & 0
\end{pmatrix},
\end{equation}
where $v_q \in \ker A_q(m^1)$ for some $q \ge p$.  It follows that 
\begin{equation} \label{eqn:upper bound}
\dim \ker A_k(m) \le \sum_{q=p}^k \dim \ker A_q(m^1).
\end{equation}

In particular, if $m_1=0$, then all elements of $\ker A_k(n)$ are of
the form $v=\begin{pmatrix}0 & \cdots & 0 & v_k\end{pmatrix}$, where
$v_k \in \ker A_k(m^1)$.  Thus,
\[
\dim \ker A_k(m)=
\dim\ker A_k(m^1)= w(m^1,|m^1|-k+1),
\]
the last equality by induction.  Since $m_1=0$, we have
$w(m,|m|-k+1)=w(m^1,|m^1|-k+1)$, which completes the proof 
in this instance.

For $0 < m_1 < k$, we consider the three cases specified in Lemma
\ref{lem:recursion}.  Assume without loss that $0 \le m_1 \le m_2 \le
\dots \le m_n$.  Let $r,r^1 \in \Z$ be maximal so that $|m|-2r \ge 0$
and $|m^1|-2r^1 \ge 0$.  Since $n \ge 3$, we have $r^1 \ge m_1$. 
Write $|m|-k+1=j$.  Since $p=k-m_1$, we have
\[
|m^1|-p+1 = j,\ 
|m^1|-(p+1)+1 = j-1,\ \dots,\ 
|m^1|-(k-1)+1=j-(m_1-1),\ 
|m^1|-k+1=j-m_1.
\]

First, consider the case $0 \le j < m_1$.  In this instance, we have
$|m^1| - q + 1 < 0$ for $p+j < q \le k$.  Theorem
\ref{thm:ker}.\ref{item:ker2} implies that $\ker A_{q}(m^1) = 0$ for
these $q$.  For $q=p+i$, $0 \le i \le p+j$, we have $|m^1|-(p+i)+1=j-i$.
Since $|m^1|-(p+j)+1=0$ and $j < m_1 \le r^1$, we have
$|m^1|=p+j-1 \ge 2r^1 \ge 2m_1 > 2j$.  So $p-1 > j$ and $p>j$.  Thus,
$0\le |m^1|-(p+i)+1<p+i$ for $0\le i \le j$.  Hence, by induction,
$\dim \ker A_{p+i}(m^1) = w(m^1,|m^1|-(p+i)+1)=w(m^1,j-i)$ for 
$0 \le i \le j$.  Consequently, in this case, \eqref{eqn:upper bound} 
yields
\[
\dim \ker A_k(m) \le \sum_{i=0}^j \dim \ker A_{p+i}(m^1) = 
\sum_{i=0}^j w(m^1,j-i) = w(m,j)
\]
by Lemma \ref{lem:recursion}.  Hence, 
$\dim \ker A_k(m)=w(m,j)=w(m,|m|-k+1)$. 
This completes the proof in the case
$0\le j < m_1$.

Next, consider the case $m_1 \le j \le r^1$.  In this instance, $|m^1|
- (p+i) + 1$ is nonnegative for each $i$, $0 \le i \le m_1$.  As
above, the condition $j \le r^1$ implies that $p>j$, and it follows
that $0\le |m^1|-(p+i)+1 < p+i$ for $i=0,1,\dots,k-p=m_1$.  So, by
induction, $\dim \ker A_{p+i}(m^1) = w(m^1,|m^1|-(p+i)+1)=w(m^1,j-i)$
for $0 \le i \le m_1$.  In this instance, \eqref{eqn:upper bound} yields
\[
\dim \ker A_k(m) \le \sum_{i=0}^{m_1}\dim \ker A_{p+i}(m^1) = 
\sum_{i=0}^{m_1} w(m^1,j-i) = w(m,j)
\]
by Lemma \ref{lem:recursion}.  Hence, 
$\dim \ker A_k(m)=w(m,j)=w(m,|m|-k+1)$. 
This completes the proof in the case
$m_1 \le j \le r^1$.

Finally, consider the case $r^1 < j \le r$.  Note that $2j-|m^1|>0$,
since $j \ge r^1$.  This, and $j=|m|-k+1$, implies that 
$|m^1|-(k-i)+1 < k-i$ for $0 \le i \le |m|-2j$, since
\[
k-i-(|m^1|-(k-i)+1)=k-j+m_1-2i \ge k-j+m_1-2|m|+4j
=k+j-|m|+2j-|m^1|=1+2j-|m^1|.
\]
By induction, $\dim \ker A_{k-i}(m^1) = w(m^1,|m^1|-(k-i)+1)$ for
$0 \le i \le |m|-2j$.  Furthermore, 
\[
\sum_{i=0}^{|m|-2j} w(m^1,|m^1|-(k-i)+1) = 
\sum_{i=0}^{|m|-2j} w(m^1,j-m_1+i)= w(m,j),
\]
by Lemma \ref{lem:recursion}.  So in this instance, 
\eqref{eqn:upper bound} yields
\[
\dim \ker A_k(m) \le \sum_{i=0}^{|m|-2j}\dim \ker A_{k-i}(m^1)\ + 
\sum_{i=|m|-2j+1}^{m_1}
\dim \ker A_{k-i}(m^1) = w(m,j)+d,
\]
where $d=\sum_{i=|m|-2j+1}^{m_1}\dim \ker A_{k-i}(m^1)$.

It may be the case that $\ker A_{k-i}(m^1)$ contains nontrivial
elements for $|m|-2j+1 \le i \le m_1$, yielding $d \neq 0$.  However,
we assert that these elements do not contribute to $\ker A_k(m)$, that
is, $\dim \ker A_k(m) = w(m,j)$.  For this, recall from
\eqref{eqn:kernel elements} that any element of $\ker A_k(m)$ is of
the form
\[
v=\begin{pmatrix} 0 & \cdots
& 0 & v_p & \cdots & v_q & 0 & \cdots & 0 \end{pmatrix},
\] 
where $v_q \in \ker A_q(m^1)$ for some $q \ge p$.  Let $K$ denote the
subspace of $\ker A_k(m)$ generated by all such vectors for which
$q=k-i$ and $0 \le i \le |m|-2j$.  Since $j=|m|-k+1$, we have 
$q \ge j+1$ for vectors $v \in K$.

Suppose that $u \in \ker A_k(m)$ and $u\notin K$.  Then,
$u=\begin{pmatrix} 0 & \cdots & 0 & u_p & \cdots & u_q & 0 & \cdots &
0 \end{pmatrix}$ with $q \le j$.  We will show that such an element
$u\in\ker A_k(m)$ is necessarily trivial.  Let $A_k^j(m)$ be the
submatrix of $A_k(m)$ given by
\[  
A^j_k(m)=
\begin{pmatrix}
D_{1,k} & A_{1}(m^1) \\
& \ddots & \ddots  \\
&& D_{p-1,k} & A_{p-1}(m^1)\\
&&&\b{0} & A_{p}(m^1)\\
&&&&D_{p+1,k} & A_{p+1}(m^1)\\
&&&&&\ddots&\ddots\\
&&&&&& D_{j,k} & A_{j}(m^1)
\end{pmatrix}.
\]
Let $\bar{u}=\begin{pmatrix} 0 & \cdots & 0 & u_p & \cdots & u_q
\end{pmatrix}$, and note that $\bar{u} \in \ker A_k^j(m)$.  As in
Corollary \ref{cor:rowcolequiv}, the matrix $A^j_k(m)$ is equivalent, 
via elementary row and column operations, to the matrix
\[
\bar{A}^j_k(m)=
\begin{pmatrix}
\b{I} & A_{1}(m^1) \\
& \ddots & \ddots  \\
&& \b{I} & A_{p-1}(m^1)\\
&&&\b{0} & A_{p}(m^1)\\
&&&& \b{I}& A_{p+1}(m^1)\\
&&&&&\ddots&\ddots\\
&&&&&& \b{I} & A_{j}(m^1)
\end{pmatrix}.
\]

Let $\bar{m}=(\bar{m}_1,\bar{m}_2,\dots,\bar{m}_n)$, where
$\bar{m}_1=j-k+m_1$ and $\bar{m}_i=m_i$ for $i \ge 2$.  Note that the
condition $j > r^1$ implies that $j-k+m_1 \ge 0$.  
By Corollary \ref{cor:rowcolequiv}, the matrix $A_j(\bar{m})$ is
equivalent to the matrix $\bar{A}^j_k(m)$ above.  Now
$|\bar{m}|-j+1=j$, so by Theorem \ref{thm:k2}.\ref{item:k23}, the
kernel of $A_j(\bar{m})$ is trivial.  Hence, the kernel of
$\bar{A}^j_k(m)$ is also trivial, and $\bar{u}=0$.  Consequently,
$u=0$, and $K=\ker A_k(m)$ as asserted.  Thus, 
$\dim \ker A_k(m) = \dim K = w(m,j)=w(m,|m|-k+1)$.  
This completes the proof in the case $r^1 < j \le r$.
\end{proof}

\end{document}